\documentclass[12pt,twoside]{article}
\usepackage{amsmath,amsfonts,amssymb,a4,enumerate,epsfig,array,theorem}

\newtheorem{theo}{Theorem}

\newtheorem{prop}{Proposition}[section]

\newtheorem{lemma}[prop]{Lemma}

\newcommand{\Eta}{E}

\newcommand{\xx}{{\bf x}}
\newcommand{\ff}{{\bf f}}
\newcommand{\OSO}{\Omega SO(3)}
\newcommand{\OSS}{\Omega \Ss^3}

\newcommand{\ZZ}{{\mathbb{Z}}}

\newcommand{\RR}{{\mathbb{R}}}

\newcommand{\Ss}{{\mathbb{S}}}
\newcommand{\BB}{{\mathbb{B}}}

\newcommand{\HH}{{\mathbb{H}}}

\newcommand{\nobf}{\noindent\bf}

\def\qed{\unskip\nobreak\hfil\penalty50\hskip1.75em\null\nobreak\hfil
$\blacksquare$ {\parfillskip=0pt \finalhyphendemerits=0 \par}\goodbreak}

\begin{document}
\title{Homotopy and cohomology of spaces \\
of locally convex curves in the sphere}
\author{Nicolau C. Saldanha}
\maketitle

\begin{abstract}
We discuss the homotopy type and the cohomology of spaces
of locally convex parametrized curves $\gamma: [0,1] \to \Ss^2$,
i.e., curves with positive geodesic curvature.
The space of all such curves with $\gamma(0) = \gamma(1) = e_1$
and $\gamma'(0) = \gamma'(1) = e_2$ is known to have three connected
components $X_{-1,c}$, $X_{1}$, $X_{-1}$.
We show several results concerning the homotopy type and cohomology
of these spaces. In particular, $X_{-1,c}$ is contractible,
$X_1$ and $X_{-1}$ are simply connected,
$\pi_2(X_{-1})$ contains a copy of $\ZZ$ and
$\pi_2(X_1)$ contains a copy of $\ZZ^2$.
Also, $H^n(X_{1}, \RR)$ and $H^n(X_{-1}, \RR)$
are nontrivial for all even $n$. More, $\dim H^{4n-2}(X_1, \RR) \ge 2$ and
$\dim H^{4n}(X_{-1}, \RR) \ge 2$ for all positive $n$.
\end{abstract}

\section{Introduction}

A curve $\gamma: [0,1] \to \Ss^2$ is called {\sl locally convex}
if its geodesic curvature is always positive, i.e., if
$\det(\gamma(t), \gamma'(t), \gamma''(t)) > 0$ for all $t$.
Let $X$ be the space of all locally convex curves
with prescribed initial point and initial direction:
$\gamma(0) = e_1$ and $\gamma'(0) = c_0 e_2$ for some $c_0 > 0$.
We define $X_Q \subset X$ to be the spaces of curves
with prescribed final point and final direction.
More precisely, for $Q \in SO(3)$, $X_Q$ is the set of $\gamma \in X$
for which $\gamma(1) = Qe_1$ and $\gamma'(1) = c_1 Qe_2$ for some $c_1 > 0$.
In particular, $X_I$ is the set of closed parametrized curves
of positive geodesic curvature with a prescribed base point
and base direction.
The topology in these spaces of curves can be taken to be $C^\infty$,
$C^k$ for some $k \ge 2$ or $H^k$ (as in Sobolev spaces) for some $k \ge 2$:
it actually makes very little difference
since it is easy to smoothen out a curve while
keeping its geodesic curvature positive.
In this paper we discuss the homotopy type and cohomology of the spaces $X_Q$.
The topology of these spaces has been discussed, among others,
by Little (\cite{Little}), B. Shapiro, M. Shapiro and Khesin
(\cite{Shapiro2}, \cite{ShapiroM}, \cite{SK})
but, as far as we could ascertain,
the main results presented here are new.

Let $Y_Q$ be the set of parametrized curves $\gamma: [0, 1] \to \Ss^2$
with $\gamma'(t) \ne 0$ for all $t$,
$\gamma(0) = e_1$, $\gamma'(0) = c_0 e_2$,
$\gamma(1) = Qe_1$, $\gamma'(1) = c_1 Qe_2$ ($c_0, c_1 > 0$).
Clearly $X_Q \subset Y_Q$.
For each $\gamma \in Y_Q$, define $\Gamma: [0,1] \to SO(3)$ by
\[ \begin{pmatrix} \gamma(t) & \gamma'(t) & \gamma''(t) \end{pmatrix} =
\Gamma(t) R(t), \]
$R(t)$ being an upper triangular matrix with positive diagonal
(the left hand side is the $3\times3$ matrix with columns
$\gamma(t)$, $\gamma'(t)$ and $\gamma''(t)$).
In other words,
\[ \Gamma(t) = \begin{pmatrix} \gamma(t) &
\widehat{\gamma'(t)} & \gamma(t) \times \widehat{\gamma'(t)}
\end{pmatrix}. \]
Recall that the universal (double) cover of $SO(3)$ is $\Ss^3 \subset \HH$,
the group of quaternions of absolute value $1$:
define $\tilde\Gamma: [0, 1] \to \Ss^3$ by
$\tilde\Gamma(0) = 1$, $\Pi \circ \tilde\Gamma = \Gamma$
where $\Pi: \Ss^3 \to SO(3)$ is the canonical projection.
This defines an injective map from $Y_Q$ to $Z_{z} \cup Z_{-z}$,
where $z$ is a quaternion with $|z| = 1$, $\Pi(z) = Q$ and
$Z_z$ is the set of continuous maps $\alpha: [0,1] \to \Ss^3$
with $\alpha(0) = 1$ and $\alpha(1) = z$.
Clearly each $Z_z$ is homotopically equivalent to $Z_1 = \OSS$,
the space of continuous maps $\alpha: [0,1] \to \Ss^3$
with $\alpha(0) = \alpha(1) = 1$.
From the Hirsch-Smale theorem (\cite{Smale}, \cite{Hirsch}),
the map $Y_Q \hookrightarrow Z_z \cup Z_{-z}$ is a homotopy equivalence.
In particular, $Y_Q$ has two connected components,
each one mapped to one of $Z_z$ and $Z_{-z}$:
call them $Y_{z}$ and $Y_{-z}$.


\begin{theo}
\label{theo:homotosur}
The inclusion $i_Q: X_Q \to Y_Q$ is homotopically surjective.
More precisely, for any compact space $K$ and any function
$f: K \to Y_Q$ there exists $g: K \to X_Q$
and a homotopy $H: [0,1] \times K \to Y_Q$
with $H(0,\cdot) = f$ and $H(1,\cdot) = g$.
\end{theo}

Theorem \ref{theo:homotosur} implies that $X_Q$
has at least two connected components, one in each of $Y_{\pm z}$.
Actually, we shall describe a set $A \subset \Ss^3$ with the following
propoerties. If $z \not\in A$, then $X \cap Y_z$ is connected:
we call this set $X_z$. If $z \in A$, then $X \cap Y_z$
has two connected components $X_z$ and $X_{z,c}$:
the inclusion $X_z \subset Y_z$ is homotopically surjective
and $X_{z,c}$ is contractible. These facts will be proved
in theorems \ref{theo:convex} and \ref{theo:XQconnect}.
Since $A$ and $-A = \{-z, z \in A\}$ will turn out to be disjoint,
$X_Q$ may have two or three components, depending on $Q$.
In \cite{Little}, Little proved that $X_I$ has three connected components.
B. Shapiro, M. Shapiro and Khesin studied the connected components
of other $X_Q$ and studied a similar problem in higher dimensions
(\cite{SK}, \cite{Shapiro2}, \cite{ShapiroM}).

We know that $H^\ast(Y_z, \RR) = H^\ast(\OSS, \RR) = \RR[\xx]$
where $\xx \in H^2(\OSS, \RR)$. Theorem \ref{theo:homotosur} implies that
the map $H^\ast(i_z): H^\ast(Y_z, \RR) \to H^\ast(X_z, \RR)$
is injective and therefore $\dim H^n(X_z, \RR) > 0$ for all even $n$.
In particular, $X_z$ is not homotopically equivalent to a finite CW-complex.

The two spaces $X_1$ and $X_{-1}$ are the objects of central interest
in this paper. We just saw that their homotopy and cohomology is at least
as large as that of $\OSS$, and we shall see that the inclusion
$X_z \subset Y_z$ is not a homotopy equivalence in either case.
On the other hand, the maps $\pi_n(i_z): \pi_n(X_z) \to \pi_n(Y_z)$
are isomorphisms for $n \le 1$.

\begin{theo}
\label{theo:simplyconnected}
The spaces $X_1$ and $X_{-1}$ are simply connected.
\end{theo}

The maps $\pi_n(i_z): \pi_n(X_z) \to \pi_n(Y_z)$
and, more generally, $[K,i_z]: [K,X_z] \to [K,Y_z]$ are surjective.
The choice of the map $g$ in theorem \ref{theo:homotosur}
is uniform up to homotopy:
intuitively, $g$ is obtained from $f$ adding many small positively
oriented loops along each $\gamma = f(k)$ so that the geodesic curvature
of $\gamma$ becomes positive (see figure \ref{fig:addloop}).
This defines maps $w_{z,K}: [K,Y_z] \to [K,X_z]$
such that $[K,i_z] \circ w_{z,K}$ is the identity.
In particular we have injective maps $w_{z,\Ss^n}: \pi_n(Y_z) \to \pi_n(X_z)$
which allow us to identify $\pi_n(Y_z)$ with a subgroup of $\pi_n(X_z)$.
Let $G_{n,z}$ be the kernel of $\pi_n(i_{z})$:
we have a natural isomorphism between $\pi_n(X_z)$
and $\pi_n(Y_z) \oplus G_{n,z}$.

For $\gamma_1 \in Y_I$ and $\gamma_2 \in Y_Q$, set
$(\gamma_1 \ast \gamma_2)(t) = \gamma_1(2t)$ for $t \le 1/2$ and
$(\gamma_1 \ast \gamma_2)(t) = \gamma_2(2t-1)$ for $t \ge 1/2$,
thus defining $\ast: Y_I \times Y_Q \to Y_Q$.
For $k > 0$, let $\nu^k \in X_I$ be given by
\[ \nu^k(t) =
\left( \frac{1 + \cos (2\pi k t)}{2}, \frac{\sqrt{2} \sin (2\pi k t)}{2},
\frac{1 - \cos (2\pi k t)}{2} \right). \]
Define $p^k_Q: X_Q \to X_Q$ by $p^k_Q(\gamma) = \nu^k \ast \gamma$.
It turns out that $p_Q = p^1_Q$ takes $X_z$ to $X_{-z}$ and vice versa.

\begin{theo}
\label{theo:homotopro}
The maps $p_Q$ and $p^3_Q$ are homotopic.
Furthermore, given $f: \Ss^n \to X_Q$ and $H: \BB^{n+1} \to Y_Q$
with $f(s) = H(s)$ for all $s \in \Ss^n$ there exists
$\tilde H: \BB^{n+1} \to X_Q$ with $\tilde H(s) = p_Q(f(s))$
for all $s \in \Ss^n$.
\end{theo}

Let $p^k_z$ be the restriction of $p^k_Q$ to $X_z \subset X_Q$.
It follows from theorem \ref{theo:homotopro} that the map
$\pi_n(p^2_z): \pi_n(X_z) \to \pi_n(X_z)$
is a projection and that the kernel of $\pi_n(p^2_z)$ is $G_{n,z}$.
We do not know what the groups $G_{n,z}$ are:
we know, however, that they are nontrivial.

\begin{theo}
\label{theo:Ghomotopy}
If $-z \in A$ then there exist functions $f_z: \Ss^2 \to X_z$
and $g_z: X_z \to \Ss^2$ such that $g_z \circ f_z$ is homotopic
to the identity. Also, $g_z \circ p^2_z$ is constant.
In particular, $G_{n,z} = \pi_n(\Ss^2) \oplus \ker(\pi_n(g_z))$.
\end{theo}

We can also say something about the cohomology of $X_z$.

\begin{theo}
\label{theo:Gcohomology}
If $(-1)^n z \in A$ then $\dim H^{2n}(X_z,\RR) \ge 2$.
\end{theo}

This work was motivated by the study of the differential equation of order 3:
\[ u'''(t) = h_1(t) u'(t) + h_0(t) u(t), \quad t \in [0,2\pi].
\eqno{(\dagger)}\]
The set of pairs of potentials $(h_0,h_1)$ for which equation $(\dagger)$
admits 3 linearly independent periodic solutions is homotopically
equivalent to $X_I$.
This appears to be the same motivation as that of
B. Shapiro and M. Shapiro for studying these same spaces.
This topological study of differential equations is
continuation of the work done together with
Dan Burghelea and Carlos Tomei in \cite{BST} and \cite{BST2}.

The author would like to thank Dan Burghelea for helpful conversations.
The author acknowledges the hospitality of
The Mathematics Department of The Ohio State University
during the winter quarter of 2004
and the support of CNPq, Capes and Faperj (Brazil).

\section{The homotopy type of $Y_Q$}

The results in this section are not new are are presented to fix notation
and for the convenience of the reader; see \cite{Arnold} for more information
concerning the geometry of curves in the sphere.

Recall that the unit tangent bundle of $\Ss^2$ is $SO(3)$:
indeed, the base point, the unit tangent vector and the cross product
of the two are the columns of an orthogonal matrix.
Also, $\pi_1(SO(3)) = \ZZ/(2)$ and the universal (double) cover of $SO(3)$
is $\Ss^3 \subset \HH$, the group of quaternions of absolute value equal to 1.
We fix notation by taking the projection $\Pi: \Ss^3 \to SO(3)$
to be given by
\[ \Pi(a + bi + cj + dk) =
\begin{pmatrix}
a^2+b^2-c^2-d^2 & -2ad+2bc & 2ac+2bd \\
2ad+2bc & a^2-b^2+c^2-d^2 & -2ab+2cd \\
-2ac+2bd & 2ab+2cd & a^2-b^2-c^2+d^2
\end{pmatrix}. \]

For an immersion $\gamma: [0,1] \to \Ss^2$,
define $\Gamma: [0, 1] \to SO(3)$ by
\[ \begin{pmatrix} \gamma(t) & \gamma'(t) & \gamma''(t) \end{pmatrix} =
\Gamma(t) R(t), \]
$R(t)$ being an upper triangular matrix with positive diagonal
(the left hand side is the $3\times3$ matrix with columns
$\gamma(t)$, $\gamma'(t)$ and $\gamma''(t)$).
In other words,
\[ \Gamma(t) = \begin{pmatrix} \gamma(t) &
\widehat{\gamma'(t)} & \gamma(t) \times \widehat{\gamma'(t)}
\end{pmatrix}. \]
If $\gamma(0) = e_1$ and $\gamma'(0) = ce_2$, $c > 0$,
define $\tilde\Gamma: [0, 1] \to \Ss^3$ by
$\tilde\Gamma(0) = 1$, $\Pi \circ \tilde\Gamma = \Gamma$.

For instance, if $\theta \in (0,\pi)$, set
\[ \gamma(t) =  \nu_\theta(t) = 
\left( \cos^2 \theta + \sin^2 \theta \cos (2\pi t),
\sin \theta \sin (2\pi t),
\cos \theta \sin \theta (1 - \cos (2\pi t)) \right). \]
The curve $\gamma$ is a circle in a plane passing through $e_1$,
parallel to $e_2$ and making an angle $\theta$ with $e_3$.
Notice that $\nu$, as defined in the introduction, is $\nu_{\pi/4}$,
that $\nu_\theta \in X_I$ for $\theta < \pi/2$ and
that $\nu_{\pi/2}$ is a geodesic. A simple computation yields
\[ \Gamma(t) =
\begin{pmatrix} \cos \theta & 0 & -\sin \theta \\
0 & 1 & 0 \\ \sin \theta & 0 & \cos \theta \end{pmatrix}
\begin{pmatrix} 1 & 0 & 0 \\ 0 & \cos(2\pi t) & -\sin(2\pi t) \\
0 & \sin(2\pi t) & \cos(2\pi t) \end{pmatrix}
\begin{pmatrix} \cos \theta & 0 & \sin \theta \\
0 & 1 & 0 \\ -\sin \theta & 0 & \cos \theta \end{pmatrix} \]
and $\tilde\Gamma(t) = \exp(\pi \ell t)$ where
$\ell = \cos(2\theta) i + \sin(2\theta) k$.
In particular, $\tilde\Gamma(1) = -1$ for all $\theta$.

The map $\phi$ taking $\gamma$ to $\Gamma$ defines a map from $Y_Q$
to $Z_Q$, the set of maps $f: [0,1] \to SO(3)$ with $f(0) = I$ and $f(1) = Q$.
Notice that $Z_Q$ is naturally identified with $Z_z \cup Z_{-z}$
where $z$ and $-z$ are the two preimages of $Q$ under $\Pi$
and $Z_z$ is the set of maps $f:[0,1] \to \Ss^3$
with $f(0) = 1$ and $f(1) = z$.
Clearly, $Z_Q$ is homotopically equivalent to $\OSO = Z_I$
and each $Z_z$ is homotopically equivalent to $\OSS = Z_1$.

Recall that $\pi_n(\OSS) = \pi_{n+1}(\Ss^3)$:
this implies that each $Z_z$ is connected
and simply connected with $\pi_2(Z_z) = \ZZ$.
Also, $H^\ast(Z_z,\RR) = H^\ast(\OSS,\RR) = \RR[\xx]$
where $\xx \in H^2(\OSS,\RR)$ satisfies $\xx^n \ne 0$ for all $n$
(see, for instance, \cite{BT}).

The Hirsch-Smale theorem proves that $\phi: Y_Q \to Z_Q$
is a homotopy equivalence: this fact admits a direct, simple
proof in our special case but we do not discuss it.
As a consequence, each $Y_z$ is connected and simply connected,
$\pi_2(Y_z) = \ZZ$ and $H^\ast(Y_z,\RR) = \RR[\xx]$.

As in the introduction, define $p^k_Q: Y_Q \to Y_Q$ by
$p^k_Q(\gamma) = \nu^k \ast \gamma$.
Notice that $p^k_Q$ as defined here
is trivially homotopic to the composition of $k$ copies of $p_Q = p^1_Q$,
justifying the notation.
The fact that $\nu \in Y_{-1}$ implies that
$p_Q: Y_Q \to Y_Q$ takes $Y_{z}$ to $Y_{-z}$ and vice-versa. 

\begin{lemma}
\label{lemma:p2i}
The function $p^2_Q: Y_{z} \to Y_{z}$ is homotopic to the identity.
\end{lemma}

{\nobf Proof: }
First notice that $p^2_Q$ is homotopic to
$\gamma \mapsto (\nu_{\epsilon} \ast \nu_{\pi - \epsilon}) \ast \gamma$
for any $\epsilon > 0$.

\begin{figure}[ht]
\begin{center}
\epsfig{height=12mm,file=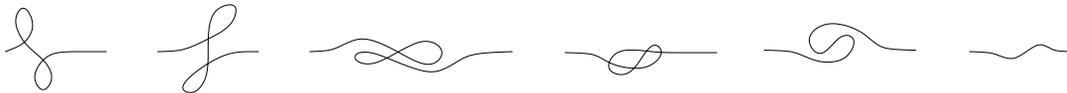}
\end{center}
\caption{How to eliminate two small loops of opposite curvature.}
\label{fig:p2i}
\end{figure}

Figure \ref{fig:p2i} shows how to move from that to a reparametrization
of $\gamma$ where a time slightly longer than $1/2$ is spent in a short
initial segment. The figure shows only the beginning of $\gamma$:
we do not think a formula is necessary or helpful.
\qed

\section{Convex curves and \\ the connected components of $X_Q$}

Given a parametrized curve of positive geodesic curvature
$\gamma: [0,1] \to \Ss^2$,
let $V_\gamma \subseteq \RR^3$ be the closure of the set
of all positive linear combinations
$a_1 \gamma(t_1) + \cdots + a_n \gamma(t_n)$,
where $a_1, \ldots, a_n$ are positive real numbers
and $t_1, \ldots, t_n \in [0, 1]$.
The set $V_\gamma$ is either a closed convex cone or $\RR^3$.
We say that $\gamma$ is {\sl convex} if
$\gamma$ is simple and the image of $\gamma$ is contained
in the boundary of $V_\gamma$.

Let $X_{Q,c} \subset X_Q$ be the set of convex curves
$\gamma: [0,1] \to \Ss^2$ with $\gamma(0) = e_1$, $\gamma'(0) = c_0 e_2$,
$\gamma(1) = Qe_1$ and $\gamma'(1) = c_1 Qe_2$ ($c_0, c_1 > 0$).
For $z$ with $\Pi(z) = Q$, define $X_z = (X_Q - X_{Q,c}) \cap Y_z$:
we shall prove in theorem \ref{theo:XQconnect} that each $X_z$ is connected.
For $Q = I$, this is proved in \cite{Little}.
In particular, $\nu \in X_{I,c}$, $\nu^{2n} \in X_1$
and $\nu^{2n+1} \in X_{-1}$ for $n$ a positive integer.


We show how to decide, given $Q$, whether $X_{Q,c}$ is empty or not.
The criterion is harder to state than to prove, so instead of proclaiming
a proposition we explain the criterion together with its justification.
We split our discussion into three cases:
$Qe_1 = e_1$, $Qe_1 = -e_1$ and $Qe_1$ and $e_1$ linearly independent.

\begin{figure}[ht]
\begin{center}
\epsfig{height=30mm,file=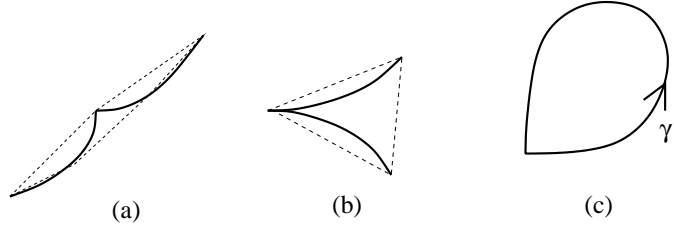}
\end{center}
\caption{How to decide if $X_{Q,c} = \emptyset$ if $Qe_1 = e_1$.}
\label{fig:XQc}
\end{figure}

If $Qe_1 = e_1$, let $\alpha \in (-\pi,\pi]$ be the angle from
$Qe_2$ to $e_2$. If $\alpha < 0$ then $X_{Q,c} = \emptyset$:
the points $\gamma(t)$ fot $t$ near $0$ or $1$ cannot possibly
be in the boundary of $V_\gamma$ if $\gamma \in X_Q$
(see figure \ref{fig:XQc}, (a): the region within the dashed line
must belong to $V_\gamma$).
Similarly, if $\alpha = \pi$ (see figure \ref{fig:XQc}, (b)).
On the other hand, if $0 \le \alpha < \pi$, it is easy to construct
$\gamma \in X_{Q,c}$ (see figure \ref{fig:XQc}. (c)).

If $Qe_1 = -e_1$, we always have $X_{Q,c} = \emptyset$.
Indeed, if a convex cone contains both $e_1$ and $-e_1$,
it must be bounded by two half-planes and there is no curve
of positive geodesic curvature contained the intersection
of these halfplanes with the unit sphere.

\begin{figure}[ht]
\begin{center}
\epsfig{height=40mm,file=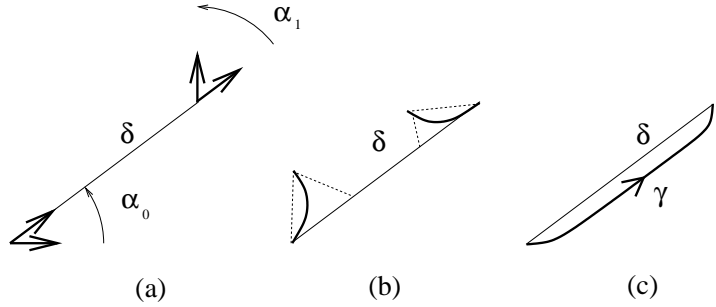}
\end{center}
\caption{How to decide if $X_{Q,c} = \emptyset$
if $e_1$ and $Qe_1$ are linearly independent.}
\label{fig:XQcdelta}
\end{figure}

Finally, consider the case when $e_1$ and $Qe_1$ are linearly independent.
Draw the shortest geodesic $\delta$ from $e_1$ to $Qe_1$:
notice that $\delta$ is contained in $V_\gamma$ for any $\gamma \in X_Q$.
Let $v_0$ and $v_1$ be the tangent vectors to $\delta$ at $e_1$ and $Qe_1$
(see figure \ref{fig:XQcdelta}).
Let $\alpha_0$ (resp. $\alpha_1$) be the angle from $e_2$ to $v_0$
(resp. from $v_1$ to $Qe_2$), $\alpha_i \in (-\pi, \pi]$.
If $\alpha_0 \le 0$ or $\alpha_1 \le 0$
then $X_{Q,c} = \emptyset$ (see figure \ref{fig:XQcdelta}, (b):
the region within the dashed line must belong to $V_\gamma$).
On the other hand, if $\alpha_0 > 0$ and $\alpha_1 > 0$
then it is easy to construct $\gamma \in X_{Q,c}$:
just keep close to $\delta$ (see figure \ref{fig:XQcdelta}, (c)).

Let $A \subset \Ss^3$ be the set of quaternions $z$ for which there
exists a convex curve in $Y_z$. The paragraphs above give a description
of $A$. The set $A$ is neither closed nor open and its interior is given by
\[ \{ a + bi + cj + dk \; | \; b, d > 0, bd > |ac| \}. \]
The sets $A$ and $-A$ are disjoint, but their closures are not.

\begin{theo}
\label{theo:convex}
Given $Q \in SO(3)$, the set $X_{Q,c} \subset X_Q$
is either empty or a contractible connected component of $X_Q$.
\end{theo}

{\nobf Proof:}
Assume that $Q$ is such that $X_{Q,c} \ne \emptyset$.
It is not hard to see that both $X_{Q,c}$ and its complement
(in $X_Q$) are open: we must prove that $X_{Q,c}$ is contractible
(and in particular that it is connected).

From the discussion above, in all cases there exists
$v \in \Ss^2$, $v \perp e_1$, $v \perp Qe_1$ such that
$\langle v, \gamma(t) \rangle > 0$ for all $\gamma \in X_{Q,c}$,
$t \in (0,1)$.
Indeed, if $Qe_1 = e_1$ we may take $v = e_3$ and
if $Qe_1 \ne e_1$ we take $v$ to be one of the vectors perpendicular
to the geodesic $\delta$.

Let $p$ be the plane $\{u | \langle v, u \rangle = 1 \}$.
For each $\gamma \in X_{Q,c}$, use radial projection
to define a curve $\hat\gamma: (0,1) \to p$
such that $\hat\gamma(t)$ is a positive multiple of $\gamma(t)$.
This defines a bijection from $X_{Q,c}$ to the set of curves
of positive curvature in the plane with prescribed asymptotic
behavior when $t$ tends to $0$ or $1$ ($\gamma(0)$ and $\gamma(1)$
indicate the asymptotic direction and $\gamma'(0)$ and $\gamma'(1)$
indicate the asymptotic line or lack thereof).
By putting axes in an appropriate position, the image of $\hat\gamma$
is the graph of a convex function from $\RR$ to $\RR$
(with prescribed asymptotic behavior).
The parametrization of the curve does not affect the homotopy type
of the space of curves, and we therefore have a homotopy equivalence
between $X_{Q,c}$ and a convex space: this proves its contractibility.
\qed

\vfil

\section{Construction of $f^{[n]}_z: (\Ss^2)^n \to X_z$}

We now give an explicit generator for $\pi_2(Y_{1})$.
This function actually has the from $f_1: \Ss^2 \to X_I \subset Y_I$
and will play an important part throughout the paper. Let
\begin{align}
\alpha_0(s,t) &= (\sin s \cos t, \sin s \sin t, \cos s), \notag\\ 
\alpha_1(s,t) &= ( -\sin t, \cos t, 0 ), \notag\\ 
\alpha_2(s,t) &= ( -\cos s \cos t, - \cos s \sin t, \sin s), \notag\\
g_s(t) &= \frac{\sqrt{2}}{2} \left( \alpha_0(s,t) + 
\cos 3t \; \alpha_1(s,t) + \sin 3t \; \alpha_2(s,t) \right). \notag
\end{align}
The curve $g_0$ is a circle drawn $4$ times and
the curve $g_\pi$ is a circle drawn $2$ times.
A computation verifies that $\det(g_s(t), g'_s(t), g''_s(t)) > 0$
for all $s$ and $t$.  Let $\Gamma_s(t)$ be defined as above.
It easy to verify that 
\[ \Gamma_s(t+(2\pi/3)) =
\begin{pmatrix} -1/2 & -\sqrt{3}/2 & 0 \\
\sqrt{3}/2 & -1/2 & 0 \\ 0 & 0 & 1 \end{pmatrix} \Gamma_s(t) \]
for all $s$ and $t$.
Finally, let $f_1: [0,2\pi] \times [0,\pi] \to X_I$ be defined by
\[ f_1(s_1, s_2)(t) =
(\Gamma_{s_2}(s_1/3))^{-1} \Gamma_{s_2}(t + (s_1/3)) e_1, \quad
s_2 \in [0,\pi] \]
If $s_2 = 0$ or $\pi$, the value of $s_1$ is irrelevant 
for the value of $f_1$:
actually, $f_1(s_1,0) = \nu^4$ and $f_1(s_1,\pi) = \nu^2$.
Also, the remark above shows that $f_1(0,s_2) = f_1(2\pi,s_2)$ for all $s_2$.
Performing these identifications, the domain of $f_1$ becomes
the sphere $\Ss^2$.

If we follow the identification of $\pi_2(Y_I)$ with $\pi_3(SO(3))$
described in the previous section, we see that in order to verify
that $f_1$ is indeed a generator of $\pi_2(Y_{I,+})$ we have to compute
the topological degree of the function
$\tilde f_1: \Ss^2 \times \Ss^1 \to \Ss^3$ which is the double cover of
$\hat f_1: \Ss^2 \times \Ss^1 \to SO(3)$ given by 
\[ \hat f_1(s,t) = \begin{pmatrix}
f_1(s)(t) & \widehat{(f_1(s))'(t)} &
f_1(s)(t) \times \widehat{(f_1(s))'(t)} \end{pmatrix}: \]
the absolute value of the degree of $\tilde f_1$ is $1$,
confirming that it is a generator.
In order to see that, it suffices to verify that $j \in \Ss^3$
is a regular value with a single preimage under $\tilde f_1$;
we skip the details.

\begin{figure}[ht]
\begin{center}
\epsfig{height=25mm,file=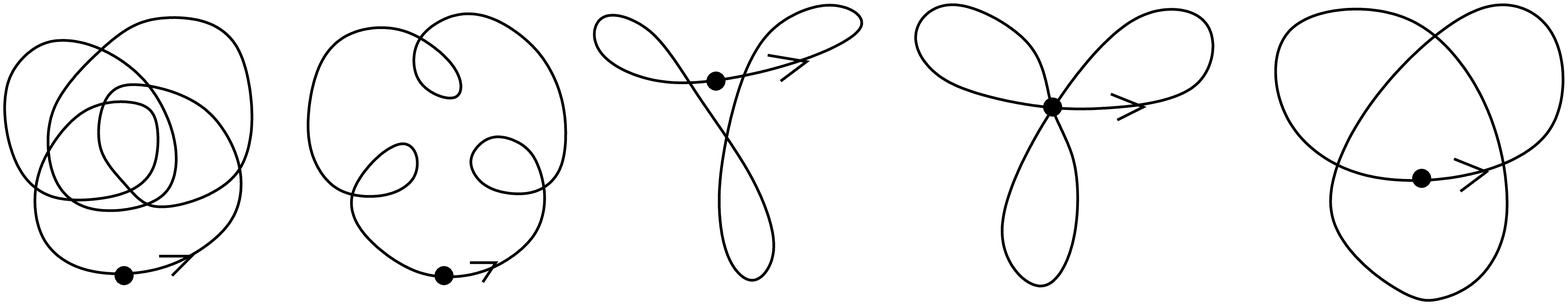}
\end{center}
\caption{The image of a line under a generator of $\pi_2(Y_{I,+})$.}
\label{fig:generator}
\end{figure}

\begin{figure}[ht]
\begin{center}
\epsfig{height=25mm,file=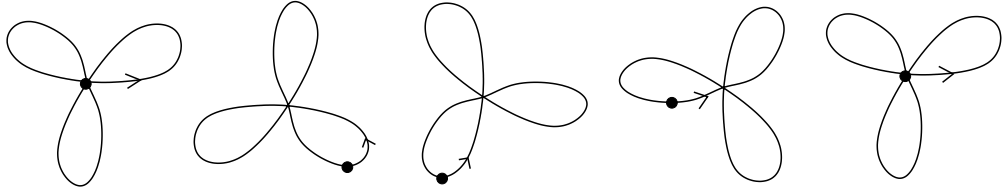}
\end{center}
\caption{The image of a circle under a generator of $\pi_2(Y_{I,+})$.}
\label{fig:otherline}
\end{figure}

Figure \ref{fig:generator} should give a rough idea of
the image of a line $(s_1,\cdot)$ under $f_1$.
The fat dot in the figure is $e_1$; between the second and third images
most of the curve went around the sphere.
Figure \ref{fig:otherline} shows the image of a line $(\cdot, s_2)$:
the first and last curves are intentionally equal:
this is a closed curve in $X_I$.

More generally, we construct, for each $z \in -A$, a function
$f_z: \Ss^2 \to X_z$ which is a generator of $\pi_2(Y_z)$.
Here, instead of a formula, we indicate the construction
in figure \ref{fig:homo}.
The dashed line is an arbitrary curve $\delta: [0,1] \to \Ss^2$
with $\delta(0) = Qe_1$ ($\Pi z = Q$) and $\delta(1) = e_1$
such that $\gamma \ast \delta$ is convex for $\gamma \in X_{-z,c}$.
The top and bottom lines are adjacent, forming a cylinder.
The curves shown are all contained in a relatively small portion
of the sphere and in the transition from the third to the fourth column,
most of the curve passed around the back of the sphere.
The way to define $f_z$ in each of the 24 squares in this grid should
be visually obvious (perhaps with a little effort).
The octagons on the right and left can likewise be filled in a natural way,
having in mind that the center of each octagon is approximately
a circle drawn two or four times.
The function $f_1$ constructed above is a special case of $f_z$.


Finally, for a positive integer $n$ and $z = (-1)^n z'$, $z' \in A$,
let $\gamma_0 \in X_{z',c}$ be fixed but arbitrary,
$\tilde\Gamma_0: [0,1] \to \Ss^3 \subset \HH$ be the associated function.
For $i = 1, 2, \ldots, n$ set
$z'_i = (\tilde\Gamma_0((i-1)/n))^{-1} \tilde\Gamma_0(i/n)$ and $z_i = - z'_i$
so that $z' = z'_1 z'_2 \cdots z'_n$ and $z = z_1 z_2 \cdots z_n$.
Notice that $z'_i \in A$. Set $Q_i = \Pi(z_i) = \Pi(z'_i)$ and define
$f^{[n]}_{z}: (\Ss^2)^n \to X_z$ by
\[ f^{[n]}_z(s_1,s_2,\ldots,s_n) =
f_{z_1}(s_1) \ast (Q_1 f_{z_2}(s_2)) \ast \cdots \ast 
(Q_1 Q_2 \cdots Q_{n-1} f_{z_n}(s_n) ). \]
Here we are stretching a bit our original definition of $\ast$:
if $\gamma_1, \gamma_2: [0,1] \to \Ss^2$
satisfy $\gamma_1(1) = \gamma_2(0)$ then 
\[ \gamma_1 \ast \gamma_2(t) = \begin{cases}
\gamma_1(2t), &t \le 1/2, \\ \gamma_2(2t-1), &t \ge 1/2. \end{cases} \]
Also, for $\gamma \in X$ and $Q \in SO(3)$,
$Q \gamma$ is the function from $[0,1]$ to $\Ss^2$ defined by
$(Q \gamma)(t) = Q (\gamma(t))$.

The following lemma is a simple consequence of the existence of $f_1$.

\begin{lemma}
\label{lemma:nu2nu4}
Let $n_1, n_2 > 1$ and $\theta_1, \theta_2 \in (0,\pi/2)$.
Then $\nu_{\theta_1}^{n_1}$ and $\nu_{\theta_2}^{n_2}$
are in the same connected component of $X_I$ if and only
if $n_1$ and $n_2$ have the same parity.
\end{lemma}

{\nobf Proof: }
First notice that if $n_1$ and $n_2$ have different parities
then $\nu_{\theta_1}^{n_1}$ and $\nu_{\theta_2}^{n_2}$
are in different connected components of $Y_I$ and therefore
with stronger reason in different connected components of $X_I$.

Our function $f_1$ shows that $\nu^2$ and $\nu^4$
are in the same connected component;
it follows from that that, for any $n > 0$,
$\nu^{2+n} \sim \nu^n \ast \nu^2$ and
$\nu^{4+n} \sim \nu^n \ast \nu^4$ are also in the same connected component.
The value of $\theta$ can be changed continuously and is therefore
not a problem. The result follows.
\qed

Recall that it follows from the results of the previous section that
$\nu \in X_{-1,c}$ is not in the same connected component of $X_I$
as $\nu^3 \in X_{-1}$.

\section{Proof of theorem \ref{theo:homotosur}}

\begin{figure}[ht]
\begin{center}
\epsfig{height=35mm,file=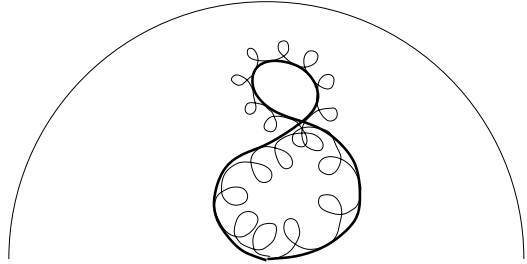}
\end{center}
\caption{A curve $\gamma$ (thicker) and $F_{9,\theta}(\gamma)$.}
\label{fig:addloop}
\end{figure}

For $\gamma \in Y_Q$ and corresponding $\Gamma: [0,1] \to SO(3)$,
define $F_{n,\theta}(\gamma)(t) = \Gamma(t) \nu_\theta^{2n}(t)$.
Intuitively, for small values of $\theta$,
$F_{n,\theta}(\gamma)$ is obtained from $\gamma$
by attaching $2n$ positively oriented small loops
along $\gamma$ (see figure \ref{fig:addloop}).

\begin{lemma}
\label{lemma:YtoX}
Let $\theta \in (0,\pi/2)$, let $K$ be a compact set and
let $f: K \to Y_Q$ a continuous function.
Then, for sufficiently large $n$,
$F_{n,\theta} \circ f$ is a function from $K$ to $X_Q$.
\end{lemma}

{\nobf Proof:}
Let $C > 1$ be a constant such that $|\Gamma'(t)| < C$
and $|\Gamma''(t)| < C$ for any $\gamma = f(k)$, $k \in K$.
Let $\epsilon > 0$ be such that if $|v_1 - \nu'_\theta(0)| < \epsilon$
and $|v_2 - \nu''_\theta(0)| < \epsilon$ then 
$\det(\nu_\theta(0), v_1, v_2) > 0$.
Notice that this implies that if $|v_1 - \nu'_\theta(t)| < \epsilon$
and $|v_2 - \nu''_\theta(t)| < \epsilon$ then 
$\det(\nu_\theta(t), v_1, v_2) > 0$.
Take $n > 20C/\epsilon$.

For $\gamma = f(k)$, write
\[ \tilde\gamma(t) = (F_{n,\theta}\gamma)(t) =
\Gamma(t) \nu^{2n}_\theta(t) = \Gamma(t) \nu_\theta(2nt) \]
so that
\begin{align}
\tilde\gamma'(t) &= \Gamma'(t) \nu_\theta(2nt) +
2n \Gamma(t) \nu_\theta'(2nt) \notag\\
\tilde\gamma''(t) &= \Gamma''(t) \nu_\theta(2nt) +
4n \Gamma'(t) \nu_\theta'(2nt) +
4n^2 \Gamma(t) \nu_\theta''(2nt) \notag
\end{align}
and therefore, after a few manipulations,
\[ \left| \frac{\tilde\gamma'(t)}{2n} - \Gamma(t)\nu_\theta'(2nt) \right|
< \epsilon, \quad
\left| \frac{\tilde\gamma''(t)}{4n^2} - \Gamma(t)\nu_\theta''(2nt) \right|
< \epsilon \]
or, equivalently,
\[ \left| \frac{(\Gamma(t))^{-1}\tilde\gamma'(t)}{2n} -
\nu_\theta'(2nt) \right| < \epsilon, \quad
\left| \frac{(\Gamma(t))^{-1}\tilde\gamma''(t)}{4n^2} -
\nu_\theta''(2nt) \right| < \epsilon. \]
It follows that
\[ \det\left(\nu_\theta(2nt),\frac{(\Gamma(t))^{-1}\tilde\gamma'(t)}{2n},
\frac{(\Gamma(t))^{-1}\tilde\gamma''(t)}{4n^2} \right) > 0 \]
and therefore that
$\det(\tilde\gamma(t), \tilde\gamma'(t), \tilde\gamma''(t)) > 0$,
which is what we needed.
\qed

Theorem \ref{theo:homotosur} now follows directly from the next lemma.

\begin{lemma}
\label{lemma:Fandp}
Let $\theta \in (0,\pi/2)$, let $K$ be a compact set, $f: K \to Y_Q$.
Then, for sufficiently large $n$,
the image of $F_{n,\theta} \circ f$ is contained in $X_Q$ and
there exists $H: [0,1] \times K \to Y_Q$ such that
$H(0,\cdot) = f$ and $H(1,\cdot) = F_{n,\theta} \circ f$.
\end{lemma}

{\nobf Proof: }
We know from lemma \ref{lemma:p2i} that $f$ if homotopic to
$f_1$, $f_1(k) = \nu_\theta^{2n} \ast f(k)$ so all we have to do
is construct a homotopy between $F_{n,\theta} \circ f$ and $f_1$.
Intuitively, this is done by pushing the loops towards $t = 0$.
More precisely, if $\gamma = f(k)$, $k \in K$, let
\[ H_1(s,k)(t) = \begin{cases} \nu_s^{2n}(t), &t \le s/2,\\
\Gamma((2t-s)/(2-s)) \nu_s^{2n}(t), &t \ge s/2 \end{cases} \]
and
\[ H_2(s,k)(t) = \begin{cases} \nu_s^{2n}((2t)/(2-s)), &t \le 1/2,\\
\Gamma(2t-1) \nu_s^{2n}((2t)/(2-s)), &1/2 \le t \le 1 - s/2,\\
\gamma(2t-1), &t \ge 1 - s/2. \end{cases} \]
Estimates similar to those of the proof of lemma \ref{lemma:YtoX}
guarantee that $H_1(s,k) \in X_Q$ and $H_2(s,k) \in Y_Q$
for sufficiently large $n$.
\qed

\section{Proof of theorem \ref{theo:homotopro}
and construction of $w_{Q,K}$}

\begin{lemma}
\label{lemma:moveloop}
Let $\theta \in (0, \pi)$.
Let $K$ be a compact space and $f: K \to Y_Q$ a continuous map.
Then, for sufficiently large $n$, if $n_1, n_2 \ge n$ then
$F_{n_1 + n_2,\theta} \circ f$ and
$\nu_{\theta}^{2n_1} \ast (F_{n_2,\theta} \circ f)$
are homotopic in the space of functions from $K$ to $X_Q$.
\end{lemma}

{\nobf Proof: }
This is a pushing-the-loops argument similar to what was done
in the proof of lemma \ref{lemma:Fandp}.
More precisely, for each $\gamma = f(k)$ first go from
$F_{n_1+n_2,\theta}(\gamma) = \Gamma \cdot \nu_\theta^{2(n_1+n_2)}$ 
to $\Gamma \cdot (\nu_\theta^{2n_1} \ast \nu_\theta^{2n_2})$:
this involves a reparametrization and the estimates in the proof
of lemma \ref{lemma:YtoX} guarantee that we remain inside $X_Q$
for sufficiently large $n_1$ and $n_2$.
Next do
\[ H(s,k)(t) = \begin{cases}
(\nu_\theta^{2n_1} \ast \nu_\theta^{2n_2})(t), &t \le s/2,\\
\Gamma((2t-s)/(2-s)) (\nu_\theta^{2n_1} \ast \nu_\theta^{2n_2})(t),
&t \ge s/2. \end{cases} \]
Again, the estimates in lemma \ref{lemma:YtoX} show that we remain
inside $X_Q$ throughout the process.
\qed

The function $w_{Q,K}: [K,Y_Q] \to [K,X_Q]$ is defined
to take $f: K \to Y_Q$ to $F_{n,\theta} \circ f: K \to X_Q$,
where $\theta$ is arbitrary and $n$ is taken to be sufficiently large.
The following lemma shows that $w_{Q,K}$ is well defined.

\begin{lemma}
\label{lemma:goodw}
Let $\theta_1, \theta_2 \in (0,\pi/2)$.
Let $K$ be a compact set and $f: K \to Y_Q$ a continuous function.
Then, there exists $N$ such that, if $n_1, n_2 > N$
then the functions $F_{n_1,\theta_1} \circ f$ and $F_{n_2, \theta_2} \circ f$
have images contained in $X_Q$ and are homotopic in the class
of functions from $K$ to $X_Q$.
\end{lemma}

{\nobf Proof: }
Take $n = \lfloor \min(n_1,n_2)/2 \rfloor$.
Use lemma \ref{lemma:moveloop} to obtain homotopies from
$F_{n_i,\theta_i} \circ f$ to
$\nu_{\theta_i}^{2(n_i-n)} \ast (F_{n,\theta_i} \circ f)$.
Lemma \ref{lemma:nu2nu4} gives us a homotopy from
$\nu_{\theta_1}^{2(n_1-n)}$ to $\nu_{\theta_2}^{2(n_2-n)}$.
A homotopy from $F_{n, \theta_1} \circ f$ to $F_{n, \theta_2} \circ f$
is obtained just by changing the value of $\theta$,
finishing the proof.
\qed

We next show that if $f: K \to X_Q$ then
$w_{Q,K} f = p^2_Q f = \nu^2 \ast f$
(the equalities here being in $[K,X_Q]$:
homotopies and not equalities as funcions).

\begin{lemma}
\label{lemma:moveloopX}
Let $\theta \in (0, \pi)$.
Let $K$ be a compact space and $f: K \to X_Q$ a continuous map.
Then, for sufficiently large $n$,
$F_{n,\theta} \circ f$ and $\nu^2 \ast f$
are homotopic in the space of functions from $K$ to $X_Q$.
\end{lemma}

{\nobf Proof: }
Write $n = n_1 + n_2$.
From lemma \ref{lemma:moveloop}, $F_{n,\theta} \circ f$
is homotopic to $\nu_\theta^{2n_1} \ast (F_{n_2,\theta} \circ f)$.
From lemma \ref{lemma:nu2nu4}, that is homotopic to
$\nu^2 \ast (F_{n_2,\theta} \circ f)$.
We now push the loops from the second interval to the first.
More precisely, for any $k \in K$, set $\gamma = \nu^2 \ast f(k)$
and corresponding $\Gamma$ so that
\[ (\nu^2 \ast (F_{n_2,\theta} \circ f))(t) = \begin{cases}
\Gamma(t) e_1, & t \le 1/2, \\
\Gamma(t) \nu_\theta^{2n_2}(2t-1), & t \ge 1/2.
\end{cases} \]
Set
\[ H(s,k)(t) = \begin{cases}
\Gamma(t) e_1, & t \le (1-s)/2, \\
\Gamma(t) \nu_\theta^{2n_2}(2t-1), & (1-s)/2 \le t \le (2-s)/2, \\
\Gamma(t) e_1, & t \ge (2-s)/2.
\end{cases} \]
This is a homotopy from $\nu^2 \ast (F_{n_2,\theta} \circ f)$
to $(F_{n_2,\theta} \circ \nu^2) \ast f$:
at any point the curvature is positive either because $n_2$
is large, using estimates similar to those of \ref{lemma:YtoX},
or because our curve is simply $\nu^2 \ast f(k)$,
which is known to be in $X_Q$.
Finally, changing the value of $\theta$ yields a homotopy
from $F_{n_2,\theta} \circ \nu^2$ to
$F_{n_2,\pi/4} \circ \nu^2 = \nu^{2(1+n_2)}$
and lemma \ref{lemma:nu2nu4} gives us a homotopy from that to $\nu^2$.
\qed

{\nobf Proof of theorem \ref{theo:homotopro}: }
The fact that $p^2_Q$ and $p^4_Q$ are homotopic follows from
lemma \ref{lemma:nu2nu4}.
Given $H: \BB^{n+1} \to Y_Q$, define $\tilde H(s) = F_{n,\theta}(H(2s))$
for $|s| \le 1/2$ where, as usual, $\theta \in (0, \pi/2)$ and $n$
is sufficiently large. Use now lemma \ref{lemma:moveloopX} for $K = \Ss^n$
to define $\tilde H$ for $1/2 < |s| < 1$.
\qed

\section{Reidemeister moves and \\ the connectivity of $X_z$}

A {\sl double point} of a curve $\gamma \in Y_I$
is a pair $(t_0, t_1)$, $0 \le t_0 < t_1 < 1$,
with $\gamma(t_0) = \gamma(t_1)$.
Similarly, a {\sl triple point} is a triple $(t_0, t_1, t_2)$,
$0 \le t_0 < t_1 < t_2 < 1$, such that
$\gamma(t_0) = \gamma(t_1) = \gamma(t_2)$.
A double point $(t_0, t_1)$ is  a {\sl self-tangency} if
$\gamma'(t_0)$ and $\gamma'(t_1)$ are parallel;
otherwise the double point is {\sl transversal}.
We call a curve {\sl generic} if it has neither triple points
nor self-tangencies and define $Y^{(0)}_I \subset Y_I$
to be the set of generic curves.
For $z = \pm 1$, set also $Y^{(0)}_z = Y^{(0)}_I \cap Y_z$,
$X^{(0)}_I = X_I \cap Y^{(0)}_I$ and $X^{(0)}_z = X^{(0)}_I \cap X_z$.
It is clear that $Y^{(0)}_z$ and $X^{(0)}_z$ are
open and dense in $Y_z$ and $X_z$, respectively.
Also, $Y^{(0)}_z$ and $X^{(0)}_z$ are disconnected
since the number of double points does not change in a connected
component of these sets.
In other words, the complement of $Y^{(0)}$ has codimension $1$.

A triple point $(t_0, t_1, t_2)$ is {\sl generic} if the vectors
$\gamma'(t_0)$, $\gamma'(t_1)$ and $\gamma'(t_2)$ are two by two
linearly independent. A self-tangency $(t_0, t_1)$ is {\sl generic}
if the curvatures of $\gamma$ at $t_0$ and $t_1$ are distinct.
A curve $\gamma \in Y_I - Y^{(0)}_I$ belongs to $Y^{(1)}_I$
if it has {\sl either} a unique generic triple point
{\sl or} a unique generic self-tangency (but not both).
The complement of $Y^{(0)}_I \cup Y^{(1)}_I \subset Y_I$ 
has codimension $2$. We define $X^{(1)}_I$, $Y^{(1)}_z$ and $X^{(1)}_z$
in the obvious way.

\begin{figure}[ht]
\begin{center}
\epsfig{height=40mm,file=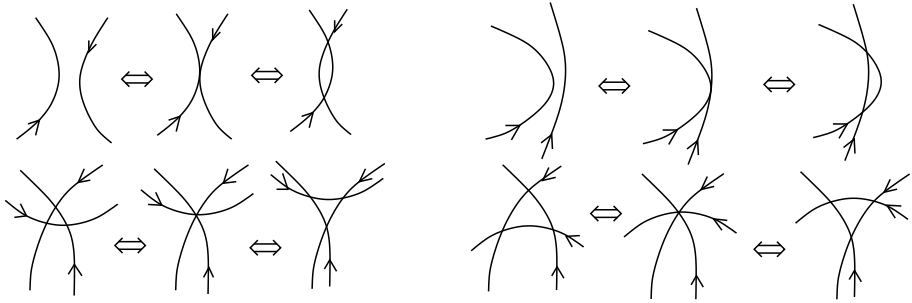}
\end{center}
\caption{Reidemeister moves: type II on first line, type III on second line.}
\label{fig:reidemeister}
\end{figure}

The passage from one connected component of $Y^{(0)}$ to another
through an element of $Y^{(1)}$ is a Reidemeister move (\cite{BZ}).
Reidemeister moves of type I are not allowed in $Y_I$;
Reidemeister moves of types II and III correspond to
generic self-tangencies and generic triple points, respectively.
If the curve is in $X_I \subset Y_I$, one of the possibilities
for orientations near a generic self-tangency is ruled out.
Figure \ref{fig:reidemeister} shows the possible Reidemeister moves,
or, equivalently, shows the neighborhood of a generic self-tangency
or generic triple point.

\begin{figure}[ht]
\begin{center}
\epsfig{height=25mm,file=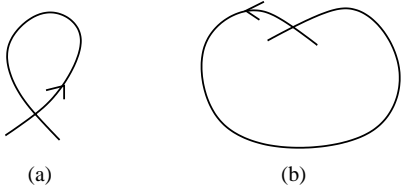}
\end{center}
\caption{A simple positive arc and a simple negative arc.}
\label{fig:littlebig}
\end{figure}

An {\it arc} of a curve $\gamma \in X_I$ is a pair $(t_-, t_+)$
such that either $(t_-,t_+)$ or $(t_+,t_-)$ is a transversal double point.
Intuitively, we think of the arc as the restriction of $\gamma$
to $[t_-,t_+]$ or $[0,t_+] \cup [t_-,1]$.
An arc is {\sl positive} (resp. {\sl negative}) if
$\det(\gamma(t_-), \gamma'(t_+), \gamma'(t_-)) > 0$ (resp. $< 0$).
An arc is {\sl simple} if one of the following two conditions holds:
\begin{enumerate}
\item{$0 \le t_- < t_+ < 1$ and
the restriction of $\gamma$ to $[t_-,t_+)$ is injective;}
\item{$0 \le t_+ < t_- < 1$ and
the restriction of $\gamma$ to $[0,t_+) \cup [t_-,1)$ is injective.}
\end{enumerate}
Figure \ref{fig:littlebig} shows examples of simple arcs.

The function $H$ constructed in the following lemma is the fundamental
building block in the proof that the spaces $X_{\pm 1}$ are connected
and simply connected.

\begin{lemma}
\label{lemma:arctoloop}
Let $(t_-, t_+)$ be a simple positive arc of $\gamma_0 \in X_I$.
Then there exists an open neighborhood $V$ of $\gamma_0$ and continuous
functions $t_-, t_+: V \to \Ss^1$ such that $(t_-(\gamma), t_+(\gamma))$
is a simple positive arc of $\gamma$ for all $\gamma \in V$.
Furthermore, there exists $H: [0,1] \times V \to X_I$
with $H(0,\gamma) = \gamma$ and $H(1,\gamma) = \nu^2 \ast \gamma$
for all $\gamma$.
\end{lemma}

\begin{figure}[ht]
\begin{center}
\epsfig{height=35mm,file=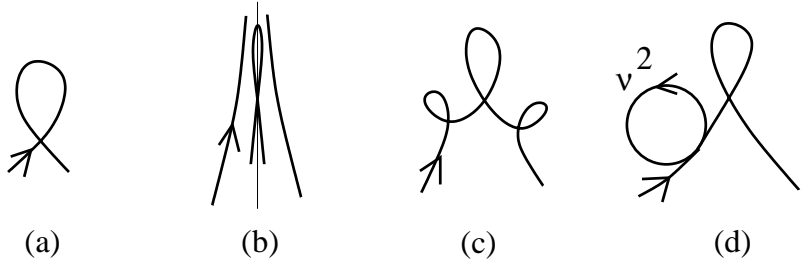}
\end{center}
\caption{How to create a copy of $\nu^2_\theta$ in a curve 
which has a small loop.}
\label{fig:little}
\end{figure}

{\nobf Proof: }
First assume $0 < t_- < t_+ < 1$.
In this case, $H(s,\gamma)$ coincides with $\gamma$ outside
$(t_- - \epsilon, t_+ + \epsilon)$ 
and $H(s,\gamma)(t)$ for $t_- - \epsilon < t < t_+ + \epsilon$
is indicated in figure \ref{fig:little}.
Let us follow the process: the arc is first shrunk (a),
then pushed along a geodesic all the way, until it comes back (b).
This creates two chunks of curve which are very nearly geodesics:
these two chunks are then shrunk (c), obtaining two new positive simple arcs
which can be deformed so that we have a copy of
$\nu^2_\theta$ (for small $\theta$) somewhere in the middle of the curve (d).
Finally, that copy of $\nu^2_\theta$ can be pushed back to $t = 0$,
proving the lemma in this case.

If $t_- = 0$ or $t_+ < t_-$, we forget about the base point
and perform the construction above to obtain
$\tilde H(s,\gamma): [0,1] \to \Ss^2$.
Define $Q(s,\gamma)$ to be the matrix in $SO(3)$ with first
two columns equal to $\tilde H(s,\gamma)(0)$ and
$\widehat{\tilde H(s,\gamma)'(0)}$ and set
$H(s,\gamma) = (Q(s,\gamma))^{-1} \tilde H(s,\gamma)$.
\qed

\begin{lemma}
\label{lemma:negarctoarc}
Let $(t_-, t_+)$ be a simple negative arc of $\gamma_0 \in X_I$.
Then there exists an open neighborhood $V$ of $\gamma_0$ and continuous
functions $t_-, t_+: V \to \Ss^1$ such that $(t_-(\gamma), t_+(\gamma))$
is a simple negative arc of $\gamma$ for all $\gamma \in V$.
Furthermore, there exists $H: [0,1] \times V \to X_I$
with $H(0,\gamma) = \gamma$ and $H(1,\gamma) = \nu^2 \ast \gamma$
for all $\gamma$.
\end{lemma}

\begin{figure}[ht]
\begin{center}
\epsfig{height=25mm,file=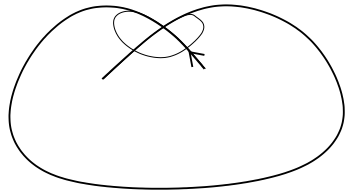}
\end{center}
\caption{How to create a simple positive arc
in a curve with a simple negative arc.}
\label{fig:big}
\end{figure}

{\nobf Proof: }
Shrink the negative simple arc as in figure \ref{fig:big}.
This creates a simple positive arc:
apply lemma \ref{lemma:arctoloop} and then unshrink.
\qed\relax

In the next result we become more general and consider
again curves in $X_z$ for $z \ne \pm 1$.

\begin{lemma}
\label{lemma:XQconnect}
Let $\gamma \in X_Q - X_{Q,c}$. Then there exists $\gamma_1 \in X_Q$
and a path in $X_Q$ from $\gamma$ to $\nu^2 \ast \gamma_1$.
\end{lemma}

{\nobf Proof: }
Assume without loss of generality that $\gamma$ is generic,
i.e., that all self intersections (if any such exist)
are transversal double (not triple) points in the interior 
(i.e., not in the endpoints). We may furthermore assume that
$\theta(t) \ne \pm e_1$ for $t \in (0,1)$.

\begin{figure}[ht]
\begin{center}
\epsfig{height=35mm,file=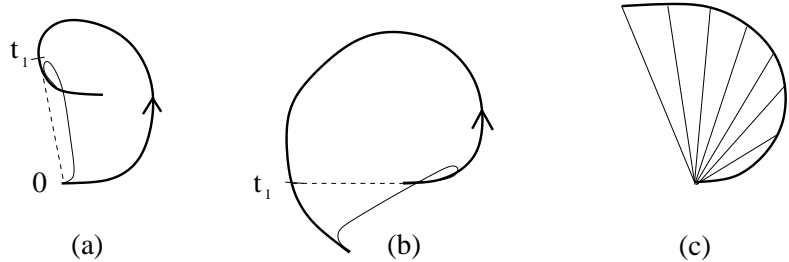}
\end{center}
\caption{How to create a self-intersection in a non-convex curve.}
\label{fig:nonconvex}
\end{figure}

If $\gamma$ is not injective, there exists a simple arc
and we use the construction in lemma \ref{lemma:arctoloop}
or \ref{lemma:negarctoarc}.
If $\gamma$ is simple but not convex, let $t_1 \in (0,1)$
be the largest number for which the restriction of $\gamma$ to $[0,t_1]$
is convex. We want to prove that this can happen in the two ways illustrated
if figure \ref{fig:nonconvex} (a) and (b). Let $\theta(t)$ be the
argument of the vector obtained from the second and third coordinates
of $\gamma(t)$: thus $\theta(t)$ is increasing for small $t$ and
its limit when $t$ tends to $0$ is $0$. As long as $\theta'(0) \ge 0$
and $\theta(t) \le \pi$, $\gamma$ restricted to $[0,t]$ is convex:
indeed, the image of this interval under $\gamma$ is a graph of
a function: one value of $x$ (the $e_1$ coordinate) for each argument
between $0$ and $\theta(t)$. In this case $V_{\gamma|_{[0,t]}} \cap \Ss^2$
is the region ``under'' this graph (see figure \ref{fig:nonconvex} (c)):
since the boundary is a locally convex simple closed curve
this set is indeed convex.

Thus $t_1$ is the first number for which either $\theta'(t) = 0$
(case (a)) or $\theta(t) = \pi$ (case (b)).
In either case it is easy to deform $\gamma$ in the interval
$[0,t_1 + \epsilon]$ as indicated in figure \ref{fig:nonconvex}.
The new (thinner) curve is very near a geodesic from $e_1$ to
$\gamma(t_1 + \epsilon)$. In either case we create self-intersections,
reducing the problem to the previous case.
\qed

\begin{theo}
\label{theo:XQconnect}
For any $z \in \Ss^3$ the set $X_z$ is connected.
\end{theo}

{\nobf Proof: }
Let $\gamma_1, \gamma_{-1} \in X_z$.
From lemma \ref{lemma:XQconnect},
each $\gamma_i$ is in the same connected component as some
$\nu^2 \ast \tilde\gamma_i$:
we must prove that
$\nu^2 \ast \tilde\gamma_1$ and $\nu^2 \ast \tilde\gamma_{-1}$
are in the same connected component of $X_z$.
Since $Y_z$ is connected, there exists
$h: \BB^1 = [-1,1] \to Y_{Q,+}$ with
$h(-1) = \tilde\gamma_{-1}$, $h(1) = \tilde\gamma_1$.
From theorem \ref{theo:homotopro}, there exists
$\tilde h: [-1, 1] \to X_Q$ with 
$\tilde h(-1) = \nu^2 \ast \tilde\gamma_{-1}$,
$\tilde h(1) = \nu^2 \ast \tilde\gamma_1$.
\qed


The construction in the proof of lemma \ref{lemma:XQconnect}
is not at all uniform.
This is not something which can be fixed with a more careful argument:
a uniform construction would prove the inclusions
$X_z \subset Y_z$ to be homotopy equivelences.
In the following sections we show that this is not the case.

\section{Proof of theorem \ref{theo:simplyconnected}}

In this section, let $z = \pm 1$.

\begin{lemma}
\label{lemma:reidhasloop}
Let $z = \pm 1$. All generic curves $\gamma \in X_z$ have simple arcs.
Furthermore, if $\gamma_0 \in X^{(1)}_z$
(a Reidemeister move) then there exists an open neighborhood
$V \subset X_z$ of $\gamma_0$ and continuous functions
$t_-, t_+: V \to \Ss^1$ such that $(t_-(\gamma), t_+(\gamma))$
is a simple arc for all $\gamma \in V$.
\end{lemma}

{\nobf Proof: }
We define a {\it generalized arc} to be a pair $(t_-,t_+)$
where $t_- \ne t_+$ and $\gamma(t_-) = \gamma(t_+)$.
We identify the generalized arc $(t_-, t_+)$ with the interval $[t_-, t_+]$
or with the set $[0,t_+] \cup [t_-,1]$,
depending on whether $t_- < t_+$ or $t_+ < t_-$ and
we order generalized arcs by inclusion.

Any curve in $\gamma \in X_z$ has self intersections.
If $\gamma$ is generic or a Reidemeister move of type III
then all generalized arcs are arcs and their number is finite;
in particular, there exists a minimal arc with respect to inclusion.
A minimal arc is clearly a simple arc.

Assume now that $\gamma$ is a Reidemeister move of type II
with self-tangency $(t_0, t_1)$. Consider the two generalized arcs
$(t_0, t_1)$ and $(t_1, t_0)$. If either is non minimal then
it contains a minimal arc and we are done. If both are simple,
they must intersect each other, otherwise we would have a Reidemeister move
from a simple closed curve in $X_{-1,c}$ to a curve in $X_{-1}$,
contradicting theorem \ref{theo:convex}. Let $(t_2, t_3)$ be this
intersection: it contains neither $(t_0, t_1)$ nor $(t_1, t_0)$.
A minimal generalized arc contained in $(t_2, t_3)$ is therefore
a simple arc.

In any of these cases the neighborhood $V$ and the functions
$t_-, t_+: V \to \Ss^1$ are constructed exactly as in
lemma \ref{lemma:arctoloop} or lemma \ref{lemma:negarctoarc}.
\qed

\begin{lemma}
\label{lemma:twoarcs}
Let $z = \pm 1$. Let $\gamma \in X_z$ be a generic curve with two
simple arcs $(t_{a,-}, t_{a,+})$ and $(t_{b,-}, t_{b,+})$.
Let $\delta_a, \delta_b: [0,1] \to X_z$ be the paths
$\delta(s) = H(s,\gamma)$ constructed in lemma \ref{lemma:arctoloop}
or \ref{lemma:negarctoarc}, so that $\delta_a(0) = \delta_b(0) = \gamma$
and $\delta_a(1) = \delta_b(1) = \nu^2 \ast \gamma$.
Then the two paths $\delta_a$ and $\delta_b$ are homotopic with
fixed endpoints in $X_z$.
\end{lemma}

{\nobf Proof: }
Notice that the statement includes the case $t_{a,\pm} = t_{b,\pm}$.
This is not quite trivial because there is an ambiguity
at the end of the construction of $H$ in lemma \ref{lemma:arctoloop}:
we did not specify which way the copy of $\nu^2$ would roll back
the base point. Thus, we have to prove that the map
$\delta_\gamma: \Ss^1 \to X_z$ taking $s_1$ to $\gamma$
with a copy of $\nu_\theta^2$ attached at the point $\gamma(s_1)$ is 
homotopic to a point. If $\theta$ is small enough, the copy of $\nu_\theta^2$
can be attached to any point of $\delta(s)$, for any $s$, thus proving
that $\delta_\gamma$ is homotopic to $\nu^2 \ast \delta_\gamma$.
Since $Y_z$ is simply connected, $\delta_\gamma$ is homotopic to a point
in $Y_z$; from theorem \ref{theo:homotopro}, $\nu^2 \ast \delta_\gamma$
is homotopic to a point in $X_z$, proving the lemma in this special case.

We next consider the case when a third simple arc $(t_{c,-}, t_{c,+})$
exists which is disjoint from the first two. Since the arcs
$(t_{a,-}, t_{a,+})$ (resp. $(t_{b,-}, t_{b,+})$) and $(t_{c,-}, t_{c,+})$
are disjoint, we may perform the construction in lemma \ref{lemma:arctoloop}
or \ref{lemma:negarctoarc} independently, thus defining
$\delta_{ac}: [0,1]^2 \to X_z$ (resp. $\delta_{bc}: [0,1]^2 to X_z$)
with $\delta_{ac}(s,0) = \delta_a(s)$ (resp. $\delta_{bc}(s,0) = \delta_b(s)$)
and $\delta_{ac}(s,1) = \nu^2 \ast \delta_a(s)$
(resp. $\delta_{bc}(s,1) = \nu^2 \ast \delta_b(s)$).
We also have $\delta_{ac}(0,s_1) = \delta_{bc}(0,s_1)$ and
$\delta_{ac}(1,s_1) = \delta_{bc}(1,s_1)$:
thus, $\delta_a$ and $\delta_b$ are homotopic with fixed endpoints
if and only if $\nu^2 \ast \delta_a$ and $\nu^2 \ast \delta_b$ are,
and this again follows from the simple connectivity of $Y_z$
and theorem \ref{theo:homotopro}.

\begin{figure}[ht]
\begin{center}
\epsfig{height=20mm,file=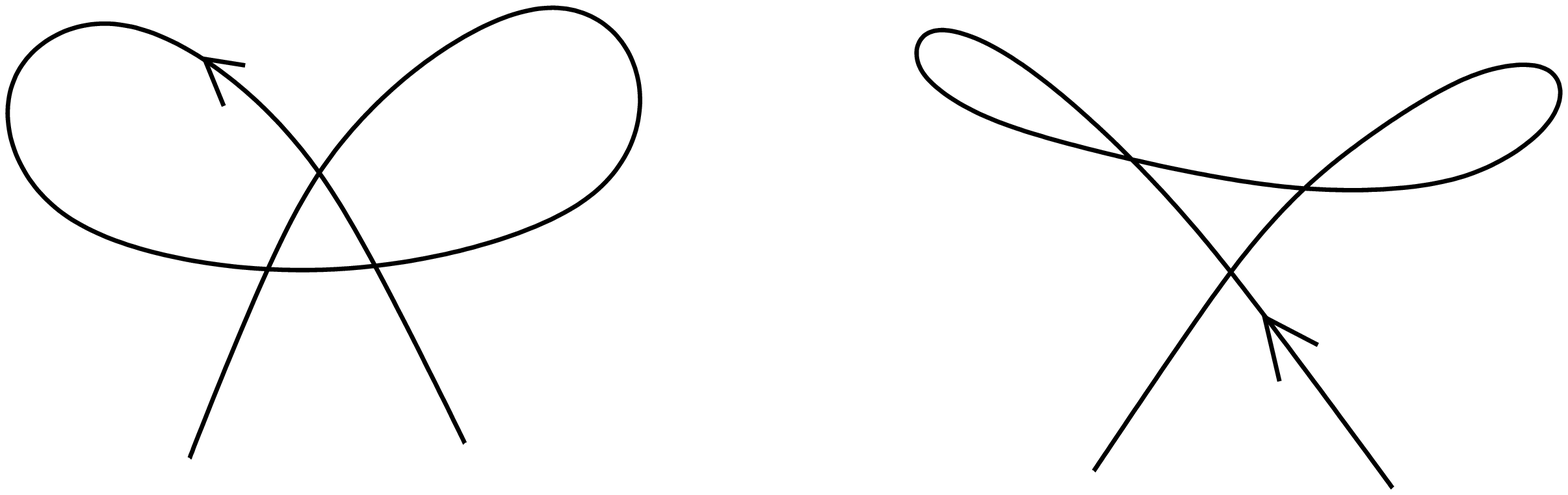}
\end{center}
\caption{Two non-disjoint positive arcs become disjoint}
\label{fig:nondisjoint}
\end{figure}

Next we consider the case when
the two arcs $(t_{a,-}, t_{a,+})$ and $(t_{b,-}, t_{b,+})$
are non-disjoint positive arcs.
In this case both $\delta_a$ and $\delta_b$ begin by performing
a Reidemeister move of type III: this makes the two arcs disjoint
and guarantees the existence of a third disjoint arc, thus reducing
the problem to the previous case, as shown in figure \ref{fig:nondisjoint}.

If one of the two initial arcs is negative, the construction
in lemma \ref{lemma:negarctoarc} creates a positive arc {\it and}
a second spare positive arc, thus again reducing the problem
to the previous cases.

\begin{figure}[ht]
\begin{center}
\epsfig{height=20mm,file=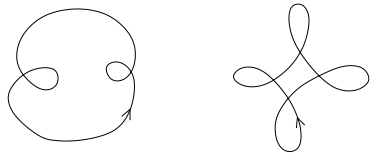}
\end{center}
\caption{If there are exactly two arcs then the argument works}
\label{fig:finalcase}
\end{figure}

We are left with only one situation to consider:
the two arcs $(t_{a,-}, t_{a,+})$ and $(t_{b,-}, t_{b,+})$
are positive, disjoint, and are the {\it only} arcs in $\gamma$.
This guarantees that, up to a deformation,
$\gamma$ is the curve in figure \ref{fig:finalcase}.
Again, both $\delta_a$ and $\delta_b$ begin by parforming a sequence
of Reidemeister moves with the same result (up to deformation),
shown in figure \ref{fig:finalcase}.
This has four disjoint positive arcs, again reducing to previous cases.
\qed

\begin{lemma}
\label{lemma:cyllinder}
Let $h: \Ss^1 \to X_z$ be a continuous function.
There exists a continuous function $H: [0,1] \times \Ss^1 \to X_z$
with $H(0,s) = h(s)$ and $H(1,s) = \nu^2 \ast (h(s))$ for all $s$.
\end{lemma}

{\nobf Proof: }
We may assume without loss of generality that $h(s)$ is generic
for all but a finite number of values of $s = s_1, \ldots, s_N$
and that these are Reidemeister moves.
Cover $\Ss^1$ by small open sets whose image under $h$ is contained
in an open neighborhood $V$ as in lemma \ref{lemma:arctoloop}
or \ref{lemma:negarctoarc}: these two results allow us to contruct
$H$ except for thin neighborhoods of finitely many transition points
from one arc to another. These transition points may be assumed
not to be Reidemeister moves. Lemma \ref{lemma:twoarcs} now guarantees
that these holes can be plugged.
\qed

{\nobf Proof of theorem \ref{theo:simplyconnected}: }
Take $h: \Ss^1 \to X_z$.
Since $Y_z$ is simply connected, $h$ is homotopic to a point in $Y_z$.
From theorem \ref{theo:homotopro}, $p_z^2 \circ h$ is homotopic
to a point in $X_z$.
From lemma \ref{lemma:cyllinder}, $h$ is homotopic to $p_z^2 \circ h$
in $X_z$. Thus, $h$ is homotopic to a point in $X_z$.
\qed

\section{Stars, trefoils and \\ the proof of theorem \ref{theo:Ghomotopy}}

\begin{figure}[ht]
\begin{center}
\epsfig{height=15mm,file=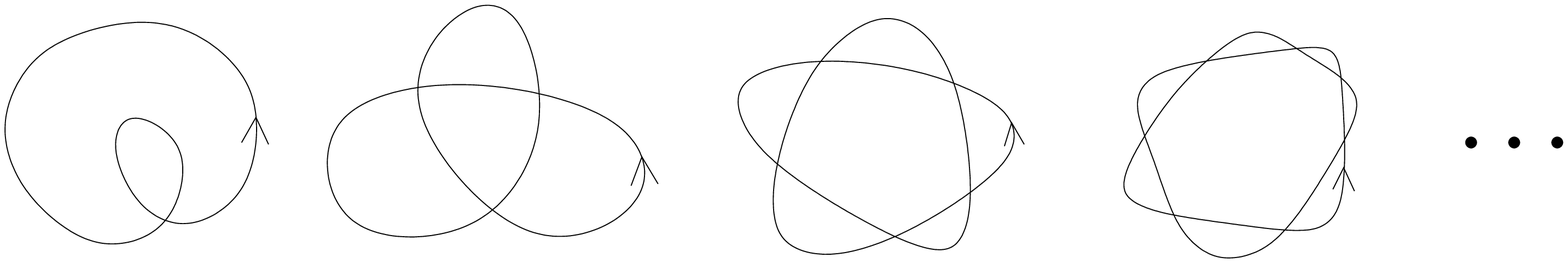}
\end{center}
\caption{Stars ($k = 0, 1, 2, 3, \ldots$).}
\label{fig:star}
\end{figure}

A {\sl star} is a curve $\gamma$ in the same connected component
of $X^{(0)}_1$ as one of the infinite family of curves given in
figure \ref{fig:star}.
More precisely, a star has $2k+1$ double points;
if $k > 0$, their images in the sphere are the vertices
of a convex polygon and, for any pair of adjacent vertices,
there are two arcs of $\gamma$ joining them.

Let $T^0$ be the closure (in $X_1$) of the set of stars
and let $T^1$ be its boundary.

A curve $\gamma \in X_1$ is called a {\sl trefoil} if:
\begin{enumerate}
\item{$\gamma$ has a generic triple point $(t_0, t_1, t_2)$;}
\item{$\gamma$ has no self-tangencies and no double points
besides $(t_0, t_1)$, $(t_0, t_2)$ and $(t_1, t_2)$;}
\item{the restriction of $\gamma$ to each of $[t_0,t_1]$, $[t_1,t_2]$
and $[t_2,1+t_0]$ is convex.}
\end{enumerate}
The restriction to $[t_2, 1 + t_0]$ is defined by $\gamma(t+1) = \gamma(t)$.
The fourth curve in figure \ref{fig:generator} and
all the curves in figure \ref{fig:otherline} are trefoils.

\begin{lemma}
\label{lemma:trefoil}
The set $T^1$ is the set of trefoils and is a manifold of codimension 1.
\end{lemma}

\begin{figure}[ht]
\begin{center}
\epsfig{height=15mm,file=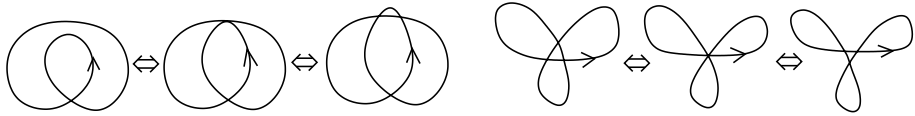}
\end{center}
\caption{Reidemeister moves starting at a star.}
\label{fig:reidstar}
\end{figure}

{\nobf Proof: }
We have to show that the only Reidemeister moves from a star
to a generic $\gamma$ which is not a star pass through a trefoil.
In order to do this, we classify all possible Reidemeister moves
starting at a star.
Figure \ref{fig:reidstar} shows how a Reidemeister move of type II
takes a star to another star (changing the value of $k$)
and how a Reidemeister move of type III takes a star ($k = 1$)
to a generic curve which is not a star passing through a trefoil.
We prove that these are the only possible moves.

\begin{figure}[ht]
\begin{center}
\epsfig{height=50mm,file=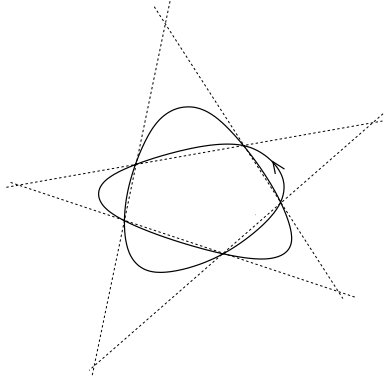}
\end{center}
\caption{A star is trapped in a union of triangles.}
\label{fig:starpoly}
\end{figure}

The only possible star from which a Reidemeister move of type III
is possible is the one shown in figure \ref{fig:reidstar} ($k = 1$):
indeed, a Reidemeister move of type III is quite impossible if the
curve does not form a combinatorial triangle.
In order to see that the only possible Reidemeister moves of type II
are those indicated in figure \ref{fig:reidstar}, notice that
if $\gamma$ is a star, its image is trapped in the union of triangles
shown in figure \ref{fig:starpoly} (where straight lines indicate
geodesics in the sphere).
\qed

Let $\tilde T^1 \subset X_1$ be an open tubular neighborhood of $T^1$
and $\Pi_{T^1}: \tilde T^1 \to T^1$ a projection onto $T^1$.
We know from lemma \ref{lemma:trefoil} that $X_1 - \tilde T^1$
has two connected components: $T^0 - \tilde T^1$ and
$X_1 - (T^0 \cup \tilde T^1)$.
Let $g_{1,a}: X_1 \to [0, \pi]$ be a continuous function 
with $g_{1,a}(\gamma) = 0$ (resp. $\pi$) for
$\gamma \in T^0 - \tilde T^1$ (resp. $X_1 - (T^0 \cup \tilde T^1)$).
We may assume without loss of generality that
$g_{1,a}(\gamma) \le \pi/2$ if and only if $\gamma \in T^0$.
Finally, let $g_1: X_1 \to \Ss^2$ be defined by
\[ g_1(\gamma) = (\sin(g_{1,a}(\gamma)) \cos\theta,
\sin(g_{1,a}(\gamma)) \sin\theta, \cos(g_{1,a}(\gamma))
\]
where, for $\gamma \in \tilde T^1$,
we have $\theta = 2\pi t_0/(1 + t_0 - t_2)$,
$(t_0, t_1, t_2)$ being the triple point of the trefoil $\Pi_{T^1}(\gamma)$.
For $\gamma \not\in \tilde T^1$, $\theta$ is undefined but this
does not affect the definition of $g_1$.

\begin{lemma}
\label{lemma:gf}
The function $g_1 \circ f_1: \Ss^2 \to \Ss^2$ has topological degree $\pm 1$.
\end{lemma}

{\nobf Proof: }
It is enough to look at the unique preimage of $(1,0,0)$,
indicated by the first and last curves in figure \ref{fig:otherline}:
the function $g_1 \circ f_1$ is injective on a neighborhood of this point.
\qed

The sign of the degree depends on the choice of orientation
for these two copies of $\Ss^2$: the domain of $f_1$ and the image of $g_1$;
we may therefore assume that these orientations were chosen so that
the degree is $1$ and then these two copies of $\Ss^2$ were identified
by an orientation preserving homeomorphism.

We wrote this section of $z = 1$ only, but everything applies
to any $z \in -A$: just attach an arc $\delta$ at the end
to make the curves closed. The definitions of star and trefoil
are based on these closed curves $\gamma \ast \delta$
and we thus define $g_z: X_z \to \Ss^2$.
The second column of figure \ref{fig:homo} consists of trefoils
and it is clear that $g_z$ takes the region to the left of this
circle to $(\cdot,\cdot,+)$, the region to the right to $(\cdot,\cdot,-)$
and that this circle is taken to the circle $(\cdot, \cdot, 0)$
by a function of degree $1$.

{\nobf Proof of theorem \ref{theo:Ghomotopy}: }
All we still have to do is prove that $g_z \circ p^2_z$ is constant:
indeed, no curve of the form $p^2_z(\gamma) = \nu^2 \ast \gamma$
will be near a star or trefoil and therefore $g_z(p^2_z(\gamma)) = (0,0,-1)$
for all $\gamma \in X_z$.
\qed

Notice that this implies that $f_z$ and $p^2_z \circ f_z$ are not homotopic.

\section{Flowers and the proof of theorem \ref{theo:Gcohomology}}

\begin{figure}[ht]
\begin{center}
\epsfig{height=30mm,file=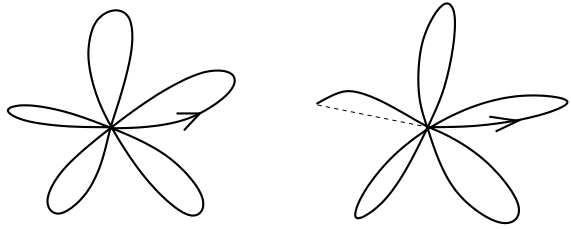}
\end{center}
\caption{Examples of flowers.}
\label{fig:flower}
\end{figure}

For $z \in A$ and $Q = \Pi(z)$, let $\theta_M \in (0,\pi]$ be
defined as follows: if $Qe_1 = e_1$, then $\theta_M$ is the
argument of $(x_2, x_3)$ where $(0, x_2, x_3) = -Qe_2$;
if $Qe_1 \ne e_1$, then $\theta_M$ is the argument of
$(x_2,x_3)$ where $Qe_1 = (x_1, x_2, x_3)$.
A curve $\gamma \in X_{(-1)^k z}$ is a {\sl flower} of $2k+1$ petals if
there exist $0 = t_0 < t_1 < t_2 < \cdots < t_{2k} < t_{2k+1} = 1$ and
$0 = \theta_0 < \theta_1 < \theta_2 < \cdots < \theta_{2k} < \theta_{2k+1}
= \theta_M$ such that:
\begin{enumerate}
\item{$\gamma(t_1) = \gamma(t_2) = \cdots = \gamma(t_{2k}) = e_1$;}
\item{the only double points of $\gamma$
are of the form $(t_i, t_j)$, $i < j$;}
\item{the argument of $(x_{i,2}, x_{i,3})$ is $\theta_i$,
where $(0, x_{i,2}, x_{i,3}) = (-1)^i \gamma'(t_i)$;}
\item{the restriction of $\gamma$ to an interval
of the form $[t_i,t_{i+1}]$ is convex.}
\end{enumerate}
Thus a flower of $1$ petal is a convex curve and
a flower of $3$ petals is a trefoil with the triple point at $e_1$.
Figure \ref{fig:flower} shows other examples of flowers.
Notice that if $\gamma$ is a flower than $\gamma(t) \ne -e_1$ for all $t$.
For $k > 0$, let $F_{2k} \subset X_{(-1)^k z}$ be the set of
flowers of $2k+1$ petals.

\begin{lemma}
\label{lemma:flower}
The set $F_{2k}$ is closed (as a subset of $X_{(-1)^k z}$
and a submanifold of codimension $2k$.
Furthermore, the normal bundle to $F_{2k}$ in $X_{(-1)^k z}$ is trivial.
\end{lemma}

\begin{figure}[ht]
\begin{center}
\epsfig{height=30mm,file=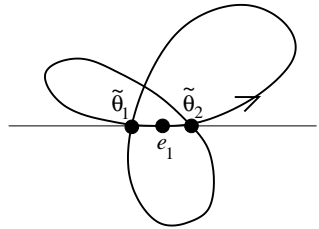}
\end{center}
\caption{A curve near a flower.}
\label{fig:almosttrefoil}
\end{figure}

{\nobf Proof: }
Any flower $\gamma$ has as open neighborhood of curves $\tilde\gamma$
as shown in figure \ref{fig:almosttrefoil}: an arc from $t = 0$
to $\tilde t_1$, with
$\tilde\gamma(\tilde t_1) = (\cos\tilde\eta_1, \sin\tilde\eta_1, 0)$,
another arc from there to $\tilde t_2$, with
$\tilde\gamma(\tilde t_2) = (\cos\tilde\eta_2, \sin\tilde\eta_2, 0)$,
and so on, 
and a final arc from $\tilde t_{2k}$ to $t_{2k+1} = 1$.
This defines a submersion $\Eta$ from this neighborhood of $\gamma$
to (a neighborhood of the origin in) $\RR^{2k}$
taking $\tilde\gamma$ to
$(\tilde\eta_1, \tilde\eta_2, \ldots, \tilde\eta_{2k})$.
Notice that $\tilde\gamma \in F_{2k+1}$ if and only if
$\Eta(\tilde\gamma) = 0$.

This construction is uniform: $\Eta$ is a submersion from a tubular
neighborhood of $F_{2k}$ to $\RR^{2k}$. This proves our claims.
\qed

The intersection number with $F_{2k}$ is therefore well defined and
can be interpreted as an element $\ff_{2k} \in H^{2k}(X_{(-1)^k z}, \ZZ)$.
We may also consider $\ff_{2k}$ to be the Poincar\'e dual of $F_{2k}$.

\begin{figure}[ht]
\begin{center}
\epsfig{height=200mm,file=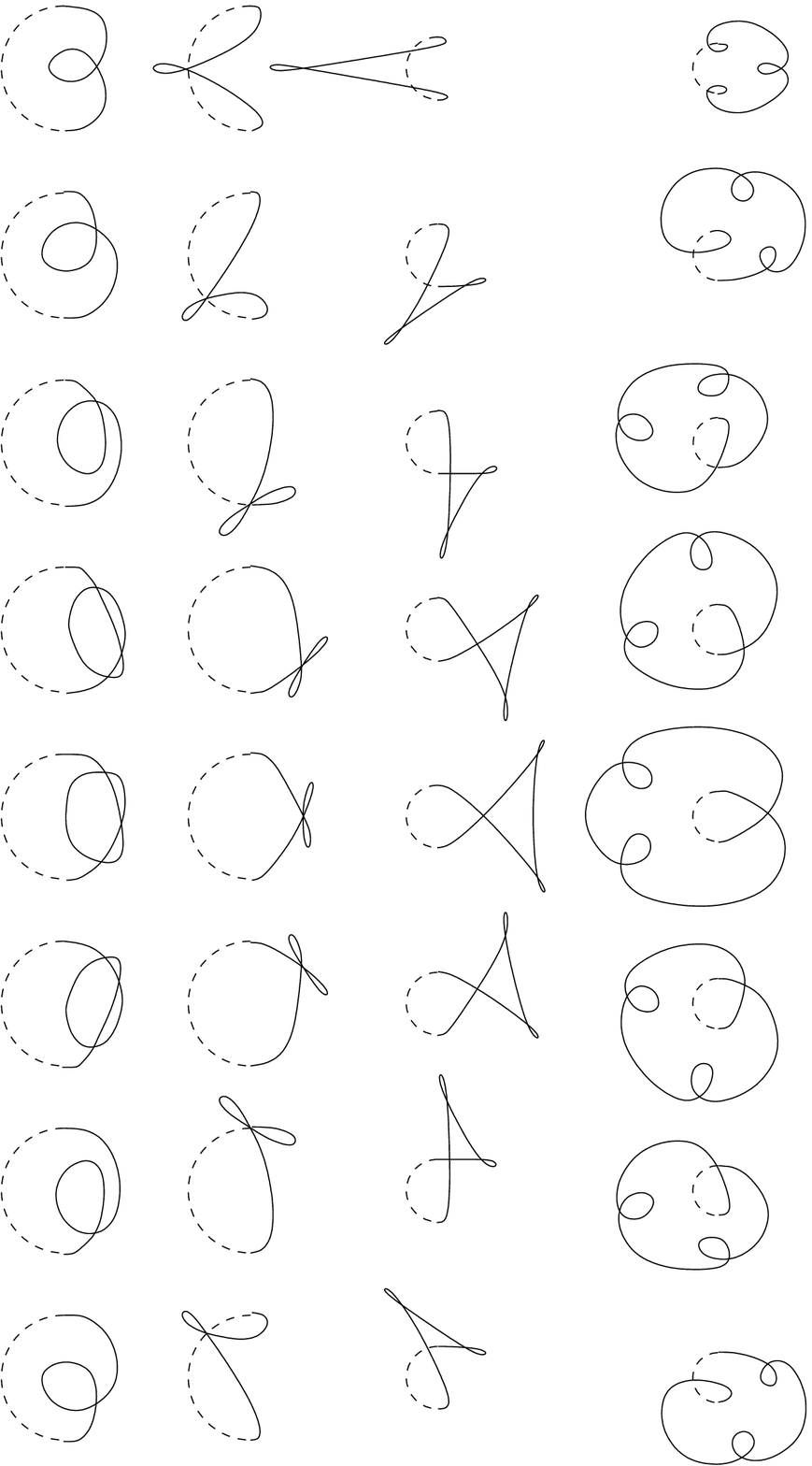}
\end{center}
\caption{The function $f_z: \Ss^2 \to X_z$, $-z \in A$.}
\label{fig:homo}
\end{figure}

{\nobf Proof of theorem \ref{theo:Gcohomology}: }
We claim that $\ff_{2k} \cdot f^{[k]}_{(-1)^k z} = \pm 1$
where $f^{[k]}_z: (\Ss^2)^k \to X_z$ was constructed in section 4
(we will again not bother with checking orientations
which we can define as we prefer anyway).
Indeed, there is a single $(s_1, s_2, \ldots , s_k) \in (\Ss^2)^k$
such that $f^{[n]}(s_1, s_2, \ldots, s_k)$ is a flower:
each $s_i \in \Ss^2$ has to be taken to be the only point
such that (in the notation of section 4) $f_{z_i}(s_i)$ is a trefoil
with the triple point at $(Q_1 Q_2 \cdots Q_{i-1})^{-1} e_1$,
which is a point on the dashed line in figure \ref{fig:homo}.
It is not hard to check that this intersection is transversal.
This proves that $\ff_{2k} \ne 0$.

On the other hand, there is obviously no flower
in the image of $\nu^2 \ast f^{[k]}_{(-1)^k z}$ and therefore
$\ff_{2k} \cdot (\nu^2 \ast f_1) = 0$.

Let $\xx \in H^2(Y_{(-1)^k z}) = H^2(\OSS)$ be as in section 2 and
let $\tilde\xx = (H^{2}(i))(\xx) \in H^{2}(X_{(-1)^k z})$.
We know from theorem \ref{theo:homotosur} that
$H^\ast(i): H^\ast(Y_{(-1)^k z}) \to H^\ast(X_{(-1)^k z})$
is injective and therefore $0 \ne \tilde\xx^k \in H^{2k}(X_{(-1)^k z})$.
By lemma \ref{lemma:p2i}, $\tilde\xx^k \cdot f^{[k]}_{(-1)^k z} =
\tilde\xx^k \cdot (\nu^2 \ast f^{[k]}_{(-1)^k z})$
and therefore $\ff_{2k}$ and $\tilde\xx^k$ are linearly independent.
\qed

We conclude by computing the product between the identified
elements of $H^\ast(X_{z},\RR)$.

\begin{prop}
\label{prop:cohomult}
For $z \in A$, let $\tilde\xx, \ff_4, \ff_8, \ldots \in H^\ast(X_z,\RR)$
be defined as above.
Then $\tilde\xx \ff_i = \ff_i \ff_j = 0$ for all $i, j$.

Similarly, for $z \in -A$, let $\tilde\xx, \ff_2, \ff_6, \ldots
\in H^\ast(X_z,\RR)$ be defined as above.
Then $\tilde\xx \ff_i = \ff_i \ff_j = 0$ for all $i, j$.
\end{prop}

{\nobf Proof: }
The product $\ff_i \ff_j$ can be interpreted in terms of the intersection
between $F_i$ and $F_j$. Since $F_i \cap F_j = \emptyset$ for $i \ne j$
it follows that $\ff_i \ff_j = 0$ in this case. Also, since the normal
bundle to $F_i$ is trivial, we can uniformly push $F_i$ in some direction
to obtain another manifold $F'_i$ homologic to $F_i$ and disjoint
from it: thus $\ff_i^2 = 0$.

In order to discuss the products $\tilde\xx \ff_i$, we recall the definition
of $\tilde\xx$. If $K$ is a closed oriented surface and $f: K \to X_z$
is a continuous function with $f(s) = \gamma$ then define
$\tilde f: K \times [0,1] \to \Ss^3$ by $\tilde f(s,t) = \tilde\Gamma(t)$.
Let $\delta: [0,1] \to \Ss^3$ be an arbitrary curve with $\delta(0) = z$,
$\delta(1) = 1$ and define $\hat f: K \times \Ss^1 \to \Ss^3$ 
(where $\Ss^1$ is $[0,1]$ with identified endpoints)
by $\hat f(s,t) = \tilde f(s,2t)$ for $t \in [0,1/2]$ and
$\hat f(s,t) = \delta(2t-1)$ for $t \in [1/2,1]$.
The product $\tilde\xx \cdot f$ is the topological degree of $\hat f$:
this can be computed using any point of $\Ss^3$ and counting its preimages
under $\hat f$ with sign. If we take that point to be $j \in \Ss^3$
and thus obtain a $2$-cocycle representing $\tilde\xx$.
Notice that having $\tilde\Gamma(t) = j$ implies $\gamma(t) = -e_1$.
Since we can not have $\gamma(t) = e_1$ for a flower,
the support of this cocycle is disjoint from $F_i$
and therefore $\tilde\xx \ff_i = 0$ for all $i$.
\qed

\vfil


\bigskip\bigskip\bigbreak

{

\parindent=0pt
\parskip=0pt
\obeylines

Nicolau C. Saldanha, PUC-Rio
nicolau@mat.puc-rio.br; http://www.mat.puc-rio.br/$\sim$nicolau/



\smallskip

Departamento de Matem\'atica, PUC-Rio
R. Marqu\^es de S. Vicente 225, Rio de Janeiro, RJ 22453-900, Brazil

}

\end{document}